\documentclass[12pt]{article}
\usepackage{amsmath}
\usepackage{amssymb}
\renewcommand\smallskip{\vskip\smallskipamount}
\renewcommand\medskip{\vskip\medskipamount}
\renewcommand\bigskip{\vskip\bigskipamount}

\begin{document}

\footnotetext{The first author acknowledges the support of NSF
Grant DMS-0654261. The second author acknowledges the support of
NSF Grant DMS-1007156 and a Sloan Research Fellowship.
}

\begin{center}
\begin{large}
\textbf{Rigidity in the Class of Orientable Compact Surfaces of
Minimal Total Absolute Curvature}
\end{large}

\bigskip\bigskip

QING HAN \text{ } \& \text{ } MARCUS KHURI

\bigskip\bigskip

\begin{abstract}
Consider an orientable compact surface in three dimensional Euclidean space with minimum total
absolute curvature. If the Gaussian curvature changes sign to finite order and satisfies a nondegeneracy condition along closed
asymptotic curves, we show that any other isometric surface differs by at most a Euclidean motion.
\end{abstract}


\textbf{1.  Introduction}
\end{center} \setcounter{equation}{0}
\setcounter{section}{1}

   A surface in $\mathbb{R}^{3}$ is said to be rigid, if any other
isometric surface differs from it by at most a Euclidean motion.
The well-known Cohn-Vossen Theorem asserts that any compact
(closed without boundary) convex surface is rigid [2].  However,
there are standard examples [8] to show that without the convexity
assumption this theorem is false.  It is then a natural question
to ask, under what conditions is a compact nonconvex surface
rigid?  Intuition suggests that a partial answer to this question
should be: a compact surface is rigid when it is as close to being
convex as is possible.  It will be our aim to verify this
assertion under a nondegeneracy condition on the Gaussian
curvature and the closed asymptotic curves.\par
   In 1938 Alexandrov [1] introduced the class $T$, of
compact surfaces (immersed in $\mathbb{R}^{3}$) characterized by
the condition
\begin{equation}
\int_{S^{+}}KdA=4\pi,
\end{equation}
where $K$ is the Gaussian curvature of $S\in T$, $S^{+}\subset S$
is the set on which $K>0$, and $dA$ is the element of area.  It is
well-known that for an arbitrary compact surface in
$\mathbb{R}^{3}$
\begin{equation*}
\int_{S}|K|dA\geq 2\pi(4-\chi(S)),
\end{equation*}
where $\chi(S)$ denotes the Euler characteristic.  The quantity on
the left-hand side is referred to as the total absolute curvature.
By the Gauss-Bonnet Theorem, if (1.1) occurs then the minimum of
the total absolute curvature is attained.  It is for this reason
that elements of the class $T$ are said to have minimal total
absolute curvature, and are sometimes referred to as tight
surfaces.  The term tightness alludes to the following geometric
consequence of (1.1),\medskip

\textbf{Convexity Property [6].}  \textit{If $S\in T$, then the
tangent plane at any point of $S^{+}$ is a plane of support, that
is, $S$ lies entirely on one side of the tangent plane.}\medskip

   Clearly, this convexity property describes the class $T$ as
those compact surfaces which are as close to being convex as is
possible for surfaces of sign changing curvature.  A standard
example of such a surface is the torus of revolution.\par
   The first rigidity result for tight surfaces was obtained by
Alexandrov [1], who showed that if $S\in T$ is analytic, then it
is rigid.  In an attempt to remove this assumption of analyticity,
Nirenberg [6] obtained the following partial result.  Let
$S^{-}\subset S\in T$ denote the set on which $K<0$.  If $\nabla
K|_{\partial S^{-}}\neq 0$, and each component of $S^{-}$ contains
at most one closed asymptotic curve, then $S$ is rigid.
Here we shall improve Nirenberg's result by
weakening the condition on the Gaussian curvature, and by allowing
multiple closed asymptotic curves with a nondegeneracy
assumption.\medskip

\textbf{Theorem 1.}  \textit{Let $S\in T$ be orientable of class
$C^{a+2}$.  If $|\nabla^{b}K|_{\partial S^{-}}\neq 0$ for some odd
integer $b<a$, and $\int_{\Gamma}k_{g}k_{n}|K|^{-1/2}ds\neq 0$ for
any closed asymptotic curve $\Gamma$ where $k_{g}$ is the geodesic
curvature, $k_{n}$ is the normal curvature in the direction
perpendicular to $\Gamma$, and $ds$ is the element of arclength,
then $S$ is rigid.}\medskip

  At the expense of adding a small extrinsic condition on
$\partial S^{-}$, we can allow the Gaussian curvature to vanish to
infinite order there.\medskip

\textbf{Theorem 2.}  \textit{Let $S\in T$ be orientable of class
$C^{4}$.  If the mean curvature $H$ satisfies $H|_{\partial
S^{-}}\neq 0$} (\textit{that is, $\partial S^{-}$ is not
umbilic}), \textit{$K$ changes sign monotonically across $\partial
S^{-}$, and $\int_{\Gamma}k_{g}k_{n}|K|^{-1/2}ds\neq 0$ for any
closed asymptotic curve $\Gamma$ where $k_{g}$ is the geodesic
curvature, $k_{n}$ is the normal curvature in the direction
perpendicular to $\Gamma$, and $ds$ is the element of arclength,
then $S$ is rigid.}\medskip

\textbf{Remark.}  \textit{At this time it is unknown whether
closed asymptotic curves in tight surfaces must always satisfy the
nondegeneracy condition of Theorems} 1 \textit{and} 2.  \textit{If
true, then Theorem} 1 \textit{would yield a purely intrinsic
condition for rigidity.}\medskip

   In order to show that a surface $S$ is rigid, it is necessary
and sufficient by the fundamental theorem of surface theory to
show that the solution of the Gauss-Codazzi equations is unique
(modulo trivial solutions).  We will use the following strategy to
verify Theorems 1 and 2.  First, a well-known argument will show
that the complement $S^{\!\text{ }'}=S-S^{-}$ is rigid.  Then
writing the Gauss-Codazzi equations as a quasi-linear weakly
hyperbolic $2\times 2$ system inside $S^{-}$ near $\partial
S^{-}$, we will use the nondegeneracy conditions placed on the
Gaussian curvature to obtain a degenerate estimate which is
sufficient to show uniqueness near $\partial S^{-}$ for the Cauchy
problem with data given on $\partial S^{-}$.  An argument of
Nirenberg [6] may then be applied to extend this local uniqueness
result inside $S^{-}$ up to a closed asymptotic curve.  Local
uniqueness for the Cauchy problem with data given on a closed
asymptotic curve will then be obtained by showing that the
Gauss-Codazzi equations form a symmetric positive system [3], from
which we can find suitable estimates.  Then repeating this
procedure throughout $S^{-}$ will yield its rigidity.
\begin{center}
\textbf{2.  The Geometry of $S^{\!\text{ }'}$ and $S^{-}$}
\end{center} \setcounter{equation}{0}
\setcounter{section}{2}

   We begin by recalling a result concerning the geometry of
$S^{\!\text{ }'}$ which will lead to its rigidity.  Let
$\mathcal{M}$ denote a 2-dimensional compact Riemannian manifold
satisfying condition (1.1), and possessing a $C^{2}$ isometric
immersion $X:\mathcal{M}\hookrightarrow\mathbb{R}^{3}$.  In [5],
Kuiper showed that $\mathcal{M}$ may be decomposed into two
disjoint open sets $U$ and $V$ with
$\mathcal{M}=\overline{U}\cup\overline{V}$ (where $\overline{U}$,
$\overline{V}$ denote the closures of $U$, $V$) such that the
restriction of $X$ to the set $U$ is an embedding and is comprised
of the boundary of the convex hull of $X(U)$ minus a finite
(possibly zero) number of planar convex disks
$D_{1},\ldots,D_{k}$.  Furthermore $K(p)>0$ for $p\in U$ and
$K(p)<0$ for $p\in V$, and the boundary of each disk $D_{i}$,
$1\leq i\leq k$, is the image of a nontrivial 1-cycle in
$\mathcal{M}$.  Now let $S_{1},S_{2}\in T$ be two isometric
immersions of $\mathcal{M}=\overline{U}\cup\overline{V}$ which
satisfy the conditions of Theorem 1 or 2, then we see that
$S_{i}^{\!\text{ }'}=X_{i}(\overline{U})$ and
$S_{i}^{-}=X_{i}(V)$, $i=1,2$, where
$X_{i}:\mathcal{M}\hookrightarrow\mathbb{R}^{3}$ are the given
immersions. Let $\gamma_{i}$ denote a boundary curve of
$S_{i}^{\!\text{ }'}$, then since $\gamma_{i}$ lies in a plane
with the normal to the surface $S_{i}$ normal to the plane, the
geodesic curvature of $\gamma_{i}$ is equal to the curvature of
$\gamma_{i}$ in the plane, which is nonnegative after appropriate
orientation.  It follows that $\gamma_{i}$ is uniquely determined
by the metric of $S_{i}$, up to a rigid body motion.  Therefore,
by filling in the empty disks in $X_{1}(\overline{U})$ and
$X_{2}(\overline{U})$ we obtain two isometric convex surfaces
$\Sigma_{1}$ and $\Sigma_{2}$.  Since $\Sigma_{1}$ and
$\Sigma_{2}$ may not be $C^{2}$ smooth, we cannot use the
Cohn-Vossen theorem to conclude that $\Sigma_{1}$ is congruent to
$\Sigma_{2}$.  However, we may apply Pogorelov's rigidity theorem
[7] for nonsmooth convex surfaces to obtain the same conclusion.
We now have
\medskip

\textbf{Lemma 2.1.}  \textit{If $S_{1},S_{2}\in T$ satisfy the
conditions of Theorem} 1 \textit{or} 2, \textit{then
$S_{1}^{\!\text{ }'}$ is congruent to $S_{2}^{\!\text{
}'}$.}\medskip

   We will now investigate the geometry of $S^{-}$ under the
hypotheses of Theorems 1 and 2.  Our goal will be to show that any
component of $S^{-}$ must be topologically equivalent to a
cylinder.  Let $\gamma$ denote one of the planar convex curves
which comprise the boundary of $S^{-}$.  Fix a point $p\in\gamma$
and introduce local coordinates $(u,v)$ near $p$, such that
$\gamma$ is the $u$-axis and the second fundamental form of $S$ is
given by
\begin{equation*}
II= Ldu^{2}+2Mdudv+Ndv^{2}.
\end{equation*}
We would like to eliminate the coefficient $M$ of the second
fundamental form.  In the case of Theorem 2, $p$ is not umbilic,
and therefore this could be done by introducing the lines of
curvature as local coordinates. However under the hypotheses of
Theorem 1, $p$ may be umbilic, so extra arguments are
needed.\medskip

\textbf{Lemma 2.2.}  \textit{Let $S$ be as in Theorem} 1
\textit{or} 2, \textit{then there exist $C^{1}$ local coordinates
$(x,t)$ near $p$, such that $\gamma$ is the $x$-axis and $M\equiv
0$ in this new coordinate system.}\medskip

\textit{Proof.}  We assume that the hypotheses of Theorem 1 hold.
Since $\gamma$ is a plane curve it has zero normal curvature and
therefore $L(u,0)=0$, this in turn implies (through the Gauss
equation (3.1)) that $M(u,0)=0$.  Then by successively
differentiating the Codazzi equations (3.1) we find that
\begin{equation}
\partial_{v}^{l}N(u,0)=0,\text{ }\text{ } l\leq k,\text{ }\text{ }
\text{ implies }\text{ }\text{ }
\partial_{v}^{l+1}L(u,0)=\partial_{v}^{l+1}M(u,0)=0,\text{ }\text{
} l\leq k.
\end{equation}
In order to obtain the desired change of coordinates set
\begin{equation*}
x=u,\text{ }\text{ }\text{ }\text{ }t=t(u,v),
\end{equation*}
where $t$ solves
\begin{equation}
t_{u}-\frac{M}{N}t_{v}=0,\text{ }\text{ }\text{ }\text{ }t(0,v)=v.
\end{equation}
By what we have just shown and by the nondegeneracy assumption on
the Gaussian curvature, the function $\frac{M}{N}$ is $C^{1}$
across the $u$-axis and satisfies $\frac{M}{N}(u,0)=0$. Therefore
the $v$-axis is noncharacteristic for (2.2), and by the theory of
first order partial differential equations it possesses a unique
$C^{1}$ local solution.  Furthermore since $t_{u}(u,0)=0$, the
$x$-axis corresponds to the curve $\gamma$. Q.E.D.\medskip

   Using Lemma 2.2 we will be able to apply an argument of
Nirenberg [6], to conclude that each component of $S^{-}$ must be
a cylinder under the hypotheses of Theorem 1 or 2.  Let $(x,t)$ be
the coordinates of Lemma 2.2 around $p\in\gamma$, and let $t>0$
denote the region lying inside $S^{-}$.  Then near $\gamma$ the
asymptotic curves of $S^{-}$, the curves which have zero normal
curvature at every point, are given by
\begin{equation*}
\frac{dt}{dx}=\pm\sqrt{-\frac{L}{N}}.
\end{equation*}
According to (2.1) $\frac{L}{N}(x,0)=0$ (under the hypotheses of
Theorem 2 this also holds, since in this case $N(x,0)\neq 0$), and
so it follows that the asymptotic curves are tangent to $\gamma$.
We will now show\medskip

\textbf{Lemma 2.3.}  \textit{Let $S$ be as in Theorem} 1
\textit{or} 2, \textit{then each component of $S^{-}$ is
topologically a cylinder.}\medskip

\textit{Proof.}  Each component of $S^{-}$ is bounded by a finite
number of planar convex curves.  Choose a component and let
$\gamma_{1},\ldots,\gamma_{k}$ denote its boundary curves. $S^{-}$
has two families of intersecting asympototic curves which are
distinguishable by how the surface rises above and below the
tangent plane at each intersection.   Take now one of these
families on the component in question.  It defines a line field
without singularity on the components closure.  By filling in the
convex boundary curves $\gamma_{1},\ldots,\gamma_{k}$ with disks
$D_{1},\ldots,D_{k}$, we obtain a closed surface on which this
line field may be extended to have a single singularity in each
disk $D_{i}$.  This extension is possible since we have shown that
the asymptotic curves are tangent to each boundary curve
$\gamma_{i}$.  The singularity in each disk clearly has index one.
Therefore, since $S$ is orientable we may apply Poincar\'{e}'s
theorem on the indices of a line field to conclude that the sum of
the indices is $2-2g$, that is $k=2-2g$ where $g$ denotes the
genus of the closed surface obtained from our component by filling
in the boundary curves with disks $D_{1},\ldots,D_{k}$.  Since
$k>0$ it follows that $g=0$ and $k=2$.  Q.E.D.\medskip

  In order to better analyze the consequences of Lemma 2.3, we
point out one consequence of the nondegeneracy condition
concerning the closed asymptotic curves in Theorems 1 and
2.\medskip

\textbf{Lemma 2.4.}  \textit{Suppose that
$\int_{\Gamma}k_{g}k_{n}|K|^{-1/2}ds\neq 0$ for a closed
asymptotic curve $\Gamma$, where $k_{g}$ is the geodesic
curvature, $k_{n}$ is the normal curvature in the direction
perpendicular to $\Gamma$, and $ds$ is the element of arclength.
Then $\Gamma$ cannot be a limit point of closed asymptotic
curves.}\medskip

\textit{Proof.}  Let $(x,t)\subset [0,x_{0})\times[0,t_{0})$ be an
arbitrary semi-global coordinate system near $\Gamma$, such that
the curve $t=0$ corresponds to $\Gamma$ and the curves
$t=\mathrm{const.}$ are all homotopic to $\Gamma$.  We will show
that
\begin{equation}
\int_{0}^{x_{0}}M^{-1}\partial_{t}L(x,0)dx=\pm\int_{\Gamma}k_{g}k_{n}|K|^{-1/2}ds\neq
0.
\end{equation}
Using the Gauss-Codazzi equations ((3.1) below) and the formulae
\begin{equation*}
k_{g}=\Gamma_{11}^{2}E^{-3/2}\sqrt{\det I},\text{ }\text{ }\text{
}\text{
}b:=\Gamma_{11}^{1}-\Gamma_{12}^{2}=2\Gamma_{11}^{1}-\frac{1}{2}\partial_{x}\log\det
I,
\end{equation*}
where $I$ denotes the first fundamental form, at $t=0$ we
calculate
\begin{eqnarray*}
M^{-1}L_{t}&=&=M^{-1}(M_{x}+bM+\Gamma_{11}^{2}N)\\
&=&b+k_{g}M^{-1}NE^{3/2}(\det I)^{-1/2}+\partial_{x}\log|M|\\
&=&2\Gamma_{11}^{1}\pm\frac{k_{g}NE^{3/2}}{\sqrt{-K}\det I}
+\partial_{x}\log|M|-\frac{1}{2}\partial_{x}\log\det I\\
&=&2\left(\Gamma_{11}^{1}\pm\frac{k_{g}FM\sqrt{E}}{\sqrt{-K}\det
I}\right)\pm\frac{k_{g}}{\sqrt{-K}}\left(\frac{EN-2FM}{\det
I}\right)\sqrt{E}\\
& &+\partial_{x}(\log|M|-\frac{1}{2}\log\det I).
\end{eqnarray*}
To discover the meaning of this expression note that
\begin{equation*}
\Gamma_{11}^{1}\pm\frac{k_{g}FM\sqrt{E}}{\sqrt{-K}\det I}
=E^{-1}(E\Gamma_{11}^{1}+F\Gamma_{11}^{2})=\frac{1}{2}\partial_{x}\log
E.
\end{equation*}
Furthermore the unit vector
\begin{equation*}
Z=-\frac{F}{\sqrt{E\det I}}\partial_{x}+\sqrt{\frac{E}{\det
I}}\partial_{t}
\end{equation*}
is normal to $\Gamma$ and satisfies
\begin{equation*}
II(Z,Z)=\frac{EN-2FM}{\det I}.
\end{equation*}
Thus
\begin{equation*}
M^{-1}L_{t}=\pm\frac{k_{g}k_{n}}{\sqrt{-K}}\sqrt{E}
+\partial_{x}\log\frac{E|M|}{\sqrt{\det I}},
\end{equation*}
so that integrating over $\Gamma$ yields (2.3).\par
  As a direct consequence of (2.3) we may confirm that $\Gamma$
cannot be a limit point of closed asymptotic curves.  Proceeding
by contradiction, assume that this is the case.  Then we may
introduce coordinates $(x,t)$, with the property that a sequence
of closed asymptotic curves $\{\Gamma_{i}\}_{i=1}^{\infty}$
converging to $\Gamma$ is given by $t=c_{i}$, for some constants
$c_{i}\rightarrow 0$ as $i\rightarrow\infty$.  It follows that
$L_{t}(x,0)=0$, in contradiction to (2.3).  Q.E.D.\medskip

   We now analyze the consequences of Lemma 2.3 with regards to
the asymptotic curves of $S^{-}$, following Nirenberg.  Let
$\gamma_{1}$ and $\gamma_{2}$ be the boundary curves for a
cylindrical component $C$ of $S^{-}$.  It is possible to orient
each family of asymptotic curves on $C$ with a suitable
parameterization, so that the corresponding vector field of
tangent vectors for each family has no singularity.  Therefore we
see that any asymptotic curve emanating from $\gamma_{1}$, say,
cannot return to $\gamma_{1}$, since otherwise a singularity in
the tangent vector field for this family must arise.  Furthermore
since $C$ is topologically a planar annulus, we may apply the
Poincar\'{e}-Bendixson Theorem to conclude that all the asymptotic
curves of one family emanating from $\gamma_{1}$ either end on
$\gamma_{2}$ or spiral towards a closed asymptotic curve
$\widetilde{\gamma}_{1}$, wrapping themselves around $C$
infinitely often.  It is clear that $\widetilde{\gamma}_{1}$ must
be homotopic to $\gamma_{1}$ and $\gamma_{2}$ since it can have no
self-intersection and cannot be homotopically trivial without
creating a singularity in the tangent vector field, and thus it
divides $C$ into two disjoint parts $C_{1}$ containing
$\gamma_{1}$ and $C_{2}$ containing $\gamma_{2}$.  Inside $C_{2}$
and near $\widetilde{\gamma}_{1}$, the curves of this family are
all spiraling towards $\widetilde{\gamma}_{1}$.  This follows from
Lemma 2.4, since $\widetilde{\gamma}_{1}$ cannot be a limit point
of closed asymptotic curves.  By applying the
Poincar\'{e}-Bendixson Theorem again, we find that these spiraling
asymptotic curves will either end on $\gamma_{2}$ or will spiral
towards another closed asymptotic curve $\widetilde{\gamma}_{2}$
homotopic to $\widetilde{\gamma}_{1}$, again dividing $C$ into two
disjoint parts $C_{1}^{\!\text{ }'}\cup C_{2}^{\!\text{ }'}$,
where $C_{2}^{\!\text{ }'}$ contains $\gamma_{2}$ and
$C_{1}^{'}=C_{1}\cup(C_{2}-C_{2}^{'})$. Inside $C_{2}^{\!\text{
}'}$ the same arguments may be applied.  This procedure may be
repeated up to $\gamma_{2}$.  Since the closure of $C$ is compact,
Lemma 2.4 implies that there can only be a finite number of closed
asymptotic curves, and so this procedure will terminate with a
finite number of iterations.  We have shown\medskip

\textbf{Lemma 2.5.}  \textit{Corresponding to each family of
asymptotic curves in $C$, is a decomposition of $C$ into a
sequence $C_{1},\ldots,C_{m}$ of cylindrical domains with
$C=\cup_{i=1}^{m}C_{i}$, such that $\partial C_{i}$, $i\geq 3$,
consists of two closed asymptotic curves of the same family
homotopic to $\gamma_{1}$ and $\gamma_{2}$, and $\partial C_{1}$,
$\partial C_{2}$ consist of one closed asymptotic curve each in
addition to $\gamma_{1}$ and $\gamma_{2}$ respectively.
Furthermore inside each $C_{i}$, $i\geq 3$, the asymptotic curves
of this family are all spiraling towards $\partial C_{i}$.  Inside
$C_{1}$, $C_{2}$ the asymptotic curves of this family are tangent
to $\gamma_{1}$, $\gamma_{2}$ and spiral towards the closed
asymptotic curve which forms the other boundary component.}

\begin{center}
\textbf{3.  The Gauss-Codazzi System}
\end{center} \setcounter{equation}{0}
\setcounter{section}{3}

   Suppose that we have two isometric surfaces
$S,\overline{S}\in T$ satisfying the assumptions of Theorem 1 or
2, and let $\mathcal{M}$ denote the underlying Riemanian manifold.
Then by Lemma 2.1, we know that $S^{\!\text{ }'}$ is congruent to
$\overline{S}^{\!\text{ }'}$.  In this section we will prove two
lemmas concerning uniqueness of the Gauss-Codazzi system which
will be instrumental in showing that $S^{-}$ is congruent to
$\overline{S}^{\text{ }\!-}$, that is we wish to show this for any
component $C$ of $S^{-}$ and the corresponding component
$\overline{C}$ for $\overline{S}^{\text{ }\!-}$.  We will denote
the corresponding cylinder in $\mathcal{M}^{-}$ by $\mathcal{C}$.
Since the boundaries of these cylinders lie in $S^{\!\text{ }'}$
and $\overline{S}^{\!\text{ }'}$, we may assume that they
coincide, after a Euclidean motion is applied.  In order to show
that $C$ is congruent to $\overline{C}$, we must show that their
second fundamental forms agree on $\mathcal{C}$.  Noting that they
are identical on the boundary, we will first prove a local (near
$\partial\mathcal{C}$) uniqueness result for the Cauchy problem of
the weakly hyperbolic system of Gauss-Codazzi equations.  This
local uniqueness will then be extended throughout $\mathcal{C}$ in
the next section.\par
   We will now put the Gauss-Codazzi system in a suitable form for
obtaining estimates.  Let $I$ and $K$ denote the metric and
Gaussian curvature of $\mathcal{C}$, and let
\begin{equation*}
II=Ldx^{2}+2Mdxdt+Ndt^{2},\text{ }\text{ }\text{ }\text{
}\overline{II}=\overline{L}dx^{2}+2\overline{M}dxdt+\overline{N}dt^{2},
\end{equation*}
denote the second fundamental forms of $C$ and $\overline{C}$ in
some local coordinate system.  Then both triples $(L,M,N)$ and
$(\overline{L},\overline{M},\overline{N})$ satisfy the
Gauss-Codazzi equations:
\begin{eqnarray}
LN-M^{2}&=&K\det I,\nonumber\\
L_{t}-M_{x}+aL+bM+cN&=&0,\\
M_{t}-N_{x}+\alpha L+\beta M+\gamma N&=&0,\nonumber
\end{eqnarray}
where $I$ denotes the first fundamental form and
$a,b,c,\alpha,\beta,\gamma$ are given in terms of Christoffel
symbols by
\begin{eqnarray*}
& &a=-\Gamma_{12}^{1},\text{ }\text{ }\text{ }\text{
}b=\Gamma_{11}^{1}-\Gamma_{12}^{2},\text{ }\text{ }\text{ }\text{
}c=\Gamma_{11}^{2},\\
& &\alpha=-\Gamma_{22}^{1},\text{ }\text{ }\text{ }\text{
}\beta=\Gamma_{12}^{1}-\Gamma_{22}^{2},\text{ }\text{ }\text{
}\text{ }\gamma=\Gamma_{12}^{2}.
\end{eqnarray*}
Set $u=\overline{L}-L$, $v=\overline{M}-M$, and
$w=\overline{N}-N$, then the triple $(u,v,w)$ satisfies the last
two (Codazzi) equations of (3.1), and in the first we may solve
for $u$ by
\begin{equation}
u=\frac{-Lw+2Mv+v^{2}}{N+w},
\end{equation}
assuming for now that the expression on the right-hand side of
(3.2) is smooth.  Plugging this into the Codazzi equations, we
have
\begin{equation*}
\left(\begin{array}{cc}
0 & -1 \\
-1 & 0 \\
\end{array}\right)\left(\begin{array}{c}
v \\
w \\
\end{array}\right)_{x}+
\left(\begin{array}{cc}
1 & 0 \\
\frac{2M}{N+w} & -\frac{L}{N+w} \\
\end{array}\right)\left(\begin{array}{c}
v \\
w \\
\end{array}\right)_{t}+
\left(\begin{array}{cc}
B_{11} & B_{12} \\
B_{21} & B_{22} \\
\end{array}\right)\left(\begin{array}{c}
v \\
w \\
\end{array}\right)=
\left(\begin{array}{c}
0 \\
0 \\
\end{array}\right),
\end{equation*}
where
\begin{eqnarray*}
B_{11}&=&\beta+\frac{\alpha v}{N+w}+\frac{2\alpha M}{N+w},\\
B_{12}&=&\gamma-\frac{\alpha L}{N+w},\\
B_{21}&=&b+\frac{2v_{t}}{N+w}-\frac{(N+w)_{t}}{(N+w)^{2}}v+2\left(\frac{M}{N+w}\right)_{t}+
\frac{2aM}{N+w}+\frac{av}{N+w},\\
B_{22}&=&c-\left(\frac{L}{N+w}\right)_{t}-\frac{aL}{N+w}.
\end{eqnarray*}
In order to symmetrize the principal portion of the system, we
multiply through by
\begin{equation*}
\left(\begin{array}{cc}
N+w & 0 \\
-2M & N+w \\
\end{array}\right)
\end{equation*}
to obtain
\begin{equation*}
\widetilde{A}^{1}U_{x}+\widetilde{A}^{2}U_{t}+\widetilde{B}U=0,
\end{equation*}
where $U=(v,w)^{*}$ (the upper $*$ will denote the transpose
operation) and
\begin{equation*}
\widetilde{A}^{1}=\left(\begin{array}{cc}
0 & -(N+w) \\
-(N+w) & 2M \\
\end{array}\right),\text{ }\text{ }\text{ }\text{ }\text{ }\text{ }\text{ }
\widetilde{A}^{2}=\left(\begin{array}{cc}
N+w & 0 \\
0 & -L \\
\end{array}\right),
\end{equation*}
\begin{equation*}
\widetilde{B}=\left(\begin{array}{cc}
(N+w)B_{11} & (N+w)B_{12} \\
-2MB_{11}+(N+w)B_{21} & -2MB_{12}+(N+w)B_{22} \\
\end{array}\right).
\end{equation*}
Lastly upon removing $w$ from the principal part the system
becomes
\begin{equation}
A^{1}U_{x}+A^{2}U_{t}+BU=0,
\end{equation}
where
\begin{equation*}
A^{1}=\left(\begin{array}{cc}
0 & -N \\
-N & 2M \\
\end{array}\right),\text{ }\text{ }\text{ }\text{ }
A^{2}=\left(\begin{array}{cc}
N & 0 \\
0 & -L \\
\end{array}\right),\text{ }\text{ }\text{ }\text{ }
B=\widetilde{B}+\left(\begin{array}{cc}
0 & v_{t}-w_{x} \\
0 & -v_{x} \\
\end{array}\right).
\end{equation*}
We will now obtain a local uniqueness result for system (3.3) near
each of the boundary curves $\Gamma_{1}$ and $\Gamma_{2}$.\medskip

\textbf{Lemma 3.1.}  \textit{There exist neighborhoods of
$\Gamma_{1}$ and $\Gamma_{2}$ inside $\mathcal{C}$ on which
$II=\overline{II}$.}\medskip

\textit{Proof.}  Let $p\in\Gamma_{1}$, and introduce the local
coordinate system $(x,t)$ of Lemma 2.2, where $p=(0,0)$ and $t>0$
represents the region inside $\mathcal{C}$.  Then $M\equiv 0$. Let
us first assume that $\partial^{k}_{t}K(x,0)\neq 0$, $k=$odd, as
in the hypothesis of Theorem 1, so that $N=O(t^{r})$ for some
$0\leq r<k$. By the arguments of Lemma 2.2, $L=O(t^{r+1})$.
Furthermore since $(\overline{L},\overline{M},\overline{N})$
agrees (to all orders) with $(L,M,N)$ on $\Gamma_{1}$, we also
have $\overline{N}=O(t^{r})$, $\overline{M}=O(t^{r+1})$, and
$\overline{L}=O(t^{r+1})$.  It follows that the expression on the
right-hand side of (3.2) is regular, so that system (3.3) is valid
near $p$.  We may assume further that $N>0$ and $L<0$ inside
$\mathcal{C}$, since $LN=K\det I$, that is
\begin{equation}
N=n(x)t^{r}+O(t^{r+1}),\text{ }\text{ }\text{ }\text{ }
L=-l(x)t^{k-r}+O(t^{k-r+1}),
\end{equation}
for some positive functions $n$ and $l$.\par
   Let $D\subset\mathcal{C}$ be a small characteristic triangle,
bounded by the $x$-axis and two intersecting asymptotic curves
given by the equations
\begin{equation}
\frac{dt}{dx}=\pm\sqrt{-\frac{L}{N}}.
\end{equation}
Note that the asymptotic curves are characteristics for the system
(3.3).  Next set $\overline{U}=e^{-\lambda}U$ where $\lambda(x,t)$
is to be determined, so that (3.3) becomes
\begin{equation}
A^{1}\overline{U}_{x}+A^{2}\overline{U}_{t}
+(B+\lambda_{x}A^{1}+\lambda_{t}A^{2})\overline{U}=0.
\end{equation}
Now multiply (3.6) by $\overline{U}^{*}$ and integrate by parts to
obtain
\begin{equation}
\int\int_{D}\overline{U}^{*}\left(\frac{\mathcal{B}+\mathcal{B}^{*}}{2}\right)\overline{U}+
\int_{\partial D}\frac{1}{2}\overline{U}^{*}(A^{1}\nu_{1}
+A^{2}\nu_{2})\overline{U}=0,
\end{equation}
where $(\nu_{1},\nu_{2})$ denotes the unit outward normal to
$\partial D$, and
\begin{equation*}
\mathcal{B}=B+\lambda_{x}A^{1}+\lambda_{t}A^{2}-\frac{1}{2}A^{1}_{x}
-\frac{1}{2}A^{2}_{t}.
\end{equation*}
In order to show that
$\frac{\mathcal{B}+\mathcal{B}^{*}}{2}:=(\mathcal{B}_{ij})$ is
positive definite, we use $v,w=O(t^{k+1})$ to calculate
\begin{eqnarray*}
\mathcal{B}_{11}&=&N\beta+\lambda_{t}N-\frac{1}{2}N_{t}+O(t^{k}),\\
\mathcal{B}_{22}&=&Nc-\frac{1}{2}L_{t}+\frac{N_{t}}{N}L-aL-\lambda_{t}L+O(t^{k}),\\
\mathcal{B}_{12}&=&\mathcal{B}_{21}=\frac{N}{2}(\gamma+b)-\frac{\alpha}{2}L
-\lambda_{x}N+\frac{1}{2}N_{x}+O(t^{k}).
\end{eqnarray*}
Moreover solving for $\gamma$ and $c$ in (3.1) yields
\begin{equation*}
\gamma=\frac{N_{x}-\alpha L}{N},\text{ }\text{ }\text{ }\text{ }
c=-\frac{L_{t}+aL}{N},
\end{equation*}
so that
\begin{equation*}
\mathcal{B}_{12}=N_{x}-\alpha
L+\frac{b}{2}N-\lambda_{x}N+O(t^{k}).
\end{equation*}
If we set $\lambda=\log
N+\frac{1}{2}\int_{0}^{x}b(\overline{x},t)d\overline{x}+\overline{\lambda}t$
for some positive constant $\overline{\lambda}$, then by (3.4)
\begin{eqnarray}
\mathcal{B}_{11}&=&N(\overline{\lambda}+\beta+\int_{0}^{x}b_{t}(\overline{x},t)d\overline{x})
+\frac{1}{2}N_{t}+O(t^{k})\nonumber\\
&=&q_{1}t^{r-1}+\overline{\lambda}q_{2}t^{r}+O(t^{r+1}),\nonumber\\
\mathcal{B}_{22}&=&-\frac{3}{2}L_{t}-(\overline{\lambda}+2a+\frac{1}{2}
\int_{0}^{x}b_{t}(\overline{x},t)d\overline{x})L+O(t^{k})\\
&=&q_{3}t^{k-r-1}+\overline{\lambda}lt^{k-r}+O(t^{k-r+1}),\nonumber\\
\mathcal{B}_{12}&=&O(t^{k-r}),\nonumber
\end{eqnarray}
where $q_{1}$, $q_{2}$, and $q_{3}$ are strictly positive
functions on $D$, and $q_{1}\equiv 0$ if $r=0$ (in which case
$\overline{\lambda}$ will be used to ensure that
$\mathcal{B}_{11}>0$).  It follows that
$\frac{\mathcal{B}+\mathcal{B}^{*}}{2}$ is positive definite since
$r<k/2$.  More precisely we have
\begin{equation}
\int\int_{D}qt^{k-r-1}|\overline{U}|^{2}\leq
\int\int_{D}\overline{U}^{*}(\frac{\mathcal{B}
+\mathcal{B}^{*}}{2})\overline{U},
\end{equation}
for some positive constant $q$.\par
   We now show that the boundary integral in (3.7) is nonnegative.
Along the $x$-axis $\overline{U}|_{t=0}=0$, since $U=O(t^{k+1})$
and $\overline{U}=e^{-\lambda}U=O(N^{-1}t^{k+1})=O(t)$.
Furthermore according to (3.5) we find that
\begin{equation*}
\nu_{1}=\pm\frac{\sqrt{-\frac{L}{N}}}{\sqrt{1-\frac{L}{N}}},\text{
}\text{ }\text{ }\text{ }\nu_{2}=\frac{1}{\sqrt{1-\frac{L}{N}}}.
\end{equation*}
Therefore
\begin{eqnarray}
& &\int_{\partial
D-\{t=0\}}\frac{1}{2}\overline{U}^{*}(A^{1}\nu_{1}
+A^{2}\nu_{2})\overline{U}\\
&=&\int_{\partial
D-\{t=0\}}\frac{1}{2}(\overline{v}\pm\sqrt{-\frac{L}{N}}\overline{w})^{2}N\nu_{2}\geq
0.\nonumber
\end{eqnarray}
By combining (3.7), (3.9), and (3.10) we conclude that $U=0$ for
$D$ sufficiently small.\par
   Now assume that the hypotheses of Theorem 2 are valid.  Then we
will slightly modify the above procedure to obtain the same
result.  We may assume that $N\geq N_{0}>0$, $L\leq 0$, and
$L_{t}|_{t>0}<0$ for $D$ small.  As above the boundary integral of
(3.7) will be nonnegative.  Moreover by choosing
$\overline{\lambda}$ sufficiently large in (3.8) we obtain
\begin{equation*}
\int\int_{D}-L_{t}|\overline{U}|^{2}\leq 0,
\end{equation*}
so that again $U=0$.  Q.E.D.\medskip

   We would now like to obtain a local uniqueness result for (3.3)
in the neighborhood of a closed asymptotic curve, $\Gamma$.  The
first step will be to construct a special semi-global coordinate
system near $\Gamma$.\medskip

\textbf{Lemma 3.2.}  \textit{Assume that
$\int_{\Gamma}k_{g}k_{n}|K|^{-1/2}ds\neq 0$ where $k_{g}$ is the
geodesic curvature, $k_{n}$ is the normal curvature in the
direction perpendicular to $\Gamma$, and $ds$ is the element of
arclength.  Then there exists a system of smooth local coordinates
$(x,t)$ near $\Gamma$ with the following properties. The
coordinate curves $t=\mathrm{const.}$ are homotopic to $\Gamma$
with $t=0$ corresponding to $\Gamma$, and the coordinate curves
$x=\mathrm{const.}$ correspond to lines of curvature all having
normal curvature of the same sign. Furthermore if}
\begin{equation*}
II=Ldx^{2}+2Mdxdt+Ndt^{2}
\end{equation*}
\textit{is the second fundamental form near $\Gamma$, then either}
\begin{equation}
L(x,0)=0,\text{ }\text{ }\text{ }\text{ }L(x,t)<0\text{ }\text{
}\textit{ for }\text{ }\text{ }t>0,\text{ }\text{ }\text{ }\text{
}\partial_{t}L(0,0)<0,
\end{equation}
\begin{equation*}
N(x,t)>0,\text{ }\text{ }\text{ }\text{ }|M(x,t)|>0,
\end{equation*}
\textit{or}
\begin{equation}
L(x,0)=0,\text{ }\text{ }\text{ }\text{ }L(x,t)>0\text{ }\text{
}\textit{ for }\text{ }\text{ }t>0,\text{ }\text{ }\text{ }\text{
}\partial_{t}L(0,0)>0,
\end{equation}
\begin{equation*}
N(x,t)<0,\text{ }\text{ }\text{ }\text{ }|M(x,t)|>0.
\end{equation*}

\textit{Proof.}  According to Lemma 2.4 $\Gamma$ cannot be a limit
point of closed asymptotic curves, and so it must be the case that
the asymptotic curves belonging to the same family as $\Gamma$ are
spiraling towards $\Gamma$.  Moreover since any asymptotic curve
of the other family (not including $\Gamma$) must intersect
$\Gamma$ transversely, and principal directions bisect asymptotic
directions, it follows that both lines of curvature intersect
$\Gamma$ transversely at every point. We may therefore choose a
preliminary semi-global coordinate system
$(\widetilde{x},\widetilde{t})$ in a sufficiently small
neighborhood of $\Gamma$ such that the curves
$\widetilde{t}=\mathrm{const.}$ are closed curves homotopic to
$\Gamma$ with $\widetilde{t}=0$ coinciding with $\Gamma$, and the
curves $\widetilde{x}=\mathrm{const.}$ are lines of curvature
corresponding to positive normal curvature.  Then the components
of the second fundamental form (in these coordinates)
$(\widetilde{L},\widetilde{M},\widetilde{N})$ satisfy
\begin{equation}
\widetilde{N}(\widetilde{x},\widetilde{t})>0,\text{ }\text{
}\text{ }\text{ }\widetilde{M}(\widetilde{x},\widetilde{t})<0,
\end{equation}
near $\Gamma$.  The estimate for $\widetilde{N}$ follows
immediately from the description of the coordinate system.
Furthermore since $\widetilde{L}(\widetilde{x},0)=0$, the Gauss
equation shows that $\widetilde{M}(\widetilde{x},0)$ cannot change
sign.  Thus (3.13) may be obtained by making the change of
coordinates $\widetilde{x}\rightarrow-\widetilde{x}$ if
necessary.\par
  We now construct the desired coordinate system.  The curves
$x=\mathrm{const}.$ remain the same, that is, they are the lines
of curvature corresponding to positive normal curvature. Choose
one of these to represent $x=0$.  The remaining curves are
labelled according to their distance along $\Gamma$ from the point
$(0,0)$. More precisely, choose an orientation for $\Gamma$ and
move along $\Gamma$ away from $(0,0)$ in the positive direction
stopping at a distance $c$. Then the (positive) line of curvature
passing through $\Gamma$ at this point will be labelled as $x=c$.
The $x$-coordinates then lie in the range $0\leq x<l$, where $l$
is the length of $\Gamma$, and all continuous functions depending
on $x$ will be periodic with period $l$.\par
  Before constructing the curves $t=\mathrm{const.}$, we make a few
observations.  In terms of the previous coordinates
$(\widetilde{x},\widetilde{t})$, the following tangent vectors
represent directions of nonpositive normal curvature:
\begin{equation*}
Y_{\sigma}=\widetilde{N}\partial_{\widetilde{x}}
+\left(-\widetilde{M}+\sigma\sqrt{-K\det\widetilde{I}}\right)
\partial_{\widetilde{t}},\text{ }\text{ }\text{ }\text{ }|\sigma|\leq 1.
\end{equation*}
In fact upon evaluating the second fundamental form in these
directions we find that
\begin{equation}
II(Y_{\sigma},Y_{\sigma})=(1-\sigma^{2})\widetilde{N}K\det\widetilde{I}.
\end{equation}
Clearly then, $Y_{\pm 1}$ are asymptotic directions and $Y_{0}$ is
a principal direction.  Furthermore since $\widetilde{M}<0$,
$Y_{-1}$ corresponds to the asymptotic direction belonging to the
same family as $\Gamma$, since this family spirals towards
$\Gamma$ and thus requires the
$\partial_{\widetilde{t}}$-component to be small.\par
  We now construct the curves $t=\mathrm{const.}$  Starting from the
point $(0,0)$, move along the curve $x=0$ (in the $t>0$ direction)
a distance $d$; this point may then be labelled $(x,t)=(0,d)$. Now
starting from $(0,d)$ follow the asymptotic curve which spirals
towards $\Gamma$ (in the positive $x$ direction), that is, the one
having direction $Y_{-1}$, until it again intersects $x=0$ at a
point $(0,d^{\!\text{ }'})$.  We assume for the time being that
$d^{\!\text{ }'}<d$, and conclude that following a curve having a
tangent vector at each point which is on the boundary of the
region with nonpositive normal curvature, produces a curve which
is not closed (since $d^{\!\text{ }'}<d$). It follows that a
smooth closed curve, starting at $(0,d)$, may be constructed by
following a curve with tangent vector at each point given by
$Y_{\sigma(x,d)}$ for appropriately chosen $\sigma(x,d)>-1$.  This
curve will be labelled by $t=d$. Analogous closed curves may be
constructed for all sufficiently small $d>0$ in such a way that
the resulting coordinate system is smooth up to $\Gamma$ ($t=0$).
Since the curves $t=\mathrm{const.}$ have tangent vector
$Y_{\sigma(x,\mathrm{const.})}$ with
$|\sigma(x,\mathrm{const.})|<1$, (3.14) shows that $L(x,t)$ is
proportional to
\begin{equation*}
(1-\sigma^{2})\widetilde{N}K\det\widetilde{I}(x,t)<0 \text{
}\text{ }\text{ for }\text{ }\text{ }t>0.
\end{equation*}
It follows that $L_{t}(x,0)\leq 0$.  Therefore (2.3) guarantees
that there exists a point on $\Gamma$ at which $L_{t}<0$, and we
can always arrange that this point occurs at $x=0$.  Lastly the
conclusion concerning $N$ is a direct consequence of the
definition of the coordinate curves, and the conclusion concerning
$M$ is a consequence of the Gauss equation and the fact that
$L(x,0)=0$.\par
  Now assume that in the construction above, $d^{\!\text{ }'}>d$.  Note that
for all sufficiently small $d$ either $d^{\!\text{ }'}>d$ or
$d^{\!\text{ }'}<d$, according to the spiraling behavior of the
asymptotic curves belonging to the same family as $\Gamma$.  In
this case we choose the curves $x=\mathrm{const.}$ to be lines of
curvature corresponding to negative normal curvature, and label
them in the same manner as described above.  To construct the
curves $t=\mathrm{const.}$, we also follow a similar procedure to
that given above.  Again, starting from the point $(x,t)=(0,d)$,
follow the asymptotic curve (having direction $Y_{-1}$) which
spirals towards $\Gamma$ in the positive $x$ direction until it
intersects $x=0$ at $(0,d^{\!\text{ }'})$. Since $d^{\!\text{
}'}>d$, we conclude that a smooth closed curve may be constructed,
starting at $(0,d)$, by following a curve with tangent vector at
each point given by $Y_{\sigma(x,d)}$ for appropriately chosen
$\sigma(x,d)<-1$.  This curve will be labelled as $t=d$. Analogous
closed curves may be constructed for all sufficiently small $d>0$
in such a way that the resulting coordinate system is smooth up to
$\Gamma$ ($t=0$). Since the curves $t=\mathrm{const.}$ now have
tangent vector $Y_{\sigma(x,\mathrm{const.})}$ with
$\sigma(x,\mathrm{const.})<-1$, (3.14) shows that $L(x,t)$ is
proportional to
\begin{equation*}
(1-\sigma^{2})\widetilde{N}K\det\widetilde{I}(x,t)>0 \text{
}\text{ }\text{ for }\text{ }\text{ }t>0.
\end{equation*}
It follows that $L_{t}(x,0)\geq 0$.  Therefore (2.3) guarantees
that there exists a point on $\Gamma$ at which $L_{t}>0$, and we
can always arrange that this point occurs at $x=0$.  Lastly the
conclusion concerning $N$ is a direct consequence of the
definition of the coordinate curves, and the conclusion concerning
$M$ is a consequence of the Gauss equation and the fact that
$L(x,0)=0$.  Q.E.D.\medskip

  We now use this special coordinate system to show that the
Gauss-Codazzi system is of symmetric positive type near the closed
asymptotic curve $\Gamma$.  This leads to the following uniqueness
result.\medskip

\textbf{Lemma 3.3.}  \textit{If $II$ and $\overline{II}$ agree on
$\Gamma$, then they agree in a neighborhood of $\Gamma$.}\medskip

\textit{Proof.}  Let $(x,t)$ be the coordinates of Lemma 3.2 and
consider the cylindrical domain $D_{\delta}=\{(x,t)\mid
0<t<\delta\}$.  We will show that $(L,M,N)$ and
$(\overline{L},\overline{M},\overline{N})$ coincide within
$D_{\delta}$ for $\delta$ sufficiently small.  In what follows we
will assume that (3.11) is valid; if (3.12) holds then nearly
identical arguments yield the desired result after multiplying
system (3.3) through by $-1$.  First observe that for sufficiently
small $\delta$ the right-hand side of (3.2) is smooth since
$N(x,t)>0$ and $w(x,0)=0$, therefore the system (3.3) is valid in
$D_{\delta}$. Set $U=e^{\lambda(x,t)+\lambda_{0}t}\overline{U}$
for some constant $\lambda_{0}>0$ and a function $\lambda(x,t)$ to
be determined, then (3.3) becomes
\begin{equation*}
A^{1}\overline{U}_{x}+A^{2}\overline{U}_{t}+\overline{B}\overline{U}=0,
\end{equation*}
where
\begin{equation*}
\overline{B}=B+\lambda_{x}A^{1}+(\lambda_{0}+\lambda_{t})A^{2}.
\end{equation*}
Multiplying by $\overline{U}^{*}$ and integrating by parts
produces
\begin{equation}
\int\int_{D_{\delta}}\overline{U}^{*}(\frac{\mathcal{B}+\mathcal{B}^{*}}{2})\overline{U}
+\int_{\partial
D_{\delta}}\frac{1}{2}\overline{U}^{*}(A^{1}\nu_{1}
+A^{2}\nu_{2})\overline{U}=0,
\end{equation}
where $(\nu_{1},\nu_{2})$ denotes the unit outer normal to
$\partial D_{\delta}$, and
\begin{equation*}
\mathcal{B}=\overline{B}-\frac{1}{2}A^{1}_{x}-\frac{1}{2}A^{2}_{t}.
\end{equation*}\par
   We now show that
$\frac{\mathcal{B}+\mathcal{B}^{*}}{2}:=(\mathcal{B}_{ij})$ is
positive definite for $\delta$ sufficiently small, if
$\lambda(x,t)$ and $\lambda_{0}$ are chosen appropriately.  A
calculation yields
\begin{eqnarray}
\mathcal{B}_{11}\!\!\!\!&=&\!\!\!\!
(\lambda_{0}+\lambda_{t}+\beta)N+\beta w-\frac{1}{2}N_{t}+\alpha
v+2\alpha M,\\
\mathcal{B}_{22}\!\!\!\!&=&\!\!\!\!
2\left(\lambda_{x}-\gamma+\frac{\alpha
L}{N+w}\right)M-M_{x}+\left(c-\left(\frac{L}{N+w}\right)_{t}\right)(N+w)\nonumber\\
& &\!\!\!\!-(\lambda_{0}+\lambda_{t}+a)L+\frac{1}{2}L_{t}-v_{x},\nonumber\\
\mathcal{B}_{12}\!=\!\mathcal{B}_{21}\!\!\!\!&=&\!\!\!\!(-\lambda_{x}
+\frac{\gamma}{2}+\frac{b}{2})N+(\frac{\gamma}{2}+\frac{b}{2})w
-\frac{\alpha}{2}L+2v_{t}-w_{x}+\frac{1}{2}(N+w)_{x}-\frac{(N+w)_{t}}{2(N+w)}v\nonumber\\
& &\!\!\!\!-\left(\beta-a+\frac{\alpha v}{N+w}+\frac{2\alpha
M}{N+w}+\frac{(N+w)_{t}}{N+w}\right)M+M_{t}+\frac{av}{2}.\nonumber
\end{eqnarray}
By applying the Codazzi equations we find that
\begin{eqnarray}
\mathcal{B}_{22}\!\!&=&\!\!(2\lambda_{x}-2\gamma)M-M_{x}+cN-\frac{1}{2}L_{t}
-(\lambda_{0}+\lambda_{t}+O(1))L+O(|w|+|\nabla w|+|\nabla
v|)\nonumber\\
&=&\!\!(2\lambda_{x}-2\gamma-b)M-\frac{3}{2}L_{t}
-(\lambda_{0}+\lambda_{t}+O(1))L+O(|w|+|\nabla w|+|\nabla
v|)\\
&=&\!\!(2\lambda_{x}-\frac{1}{2}\partial_{x}\log\det
I)M-\frac{3}{2}L_{t}
-(\lambda_{0}+\lambda_{t}+O(1))L+O(|w|+|\nabla w|+|\nabla
v|),\nonumber
\end{eqnarray}
where we have also used the identity
\begin{equation*}
2\gamma+b=\Gamma_{11}^{1}+\Gamma_{12}^{2}=\frac{1}{2}\partial_{x}\log\det
I.
\end{equation*}
This motivates the choice
\begin{equation*}
\lambda:=\lambda_{1}+\varepsilon\lambda_{2},\text{ }\text{ }\text{
}\text{ } \lambda_{1}=\frac{1}{4}\log\det I,
\end{equation*}
where $\varepsilon>0$ is a small parameter and $\lambda_{2}$ is
required to satisfy $\partial_{x}\lambda_{2}>0$ if $M(x,t)>0$
($\partial_{x}\lambda_{2}<0$ if $M(x,t)<0$) except in a
sufficiently small neighborhood of $(0,0)$ where $L_{t}<0$. Note
that since all functions, and in particular $\lambda_{2}$, must be
periodic in $x$, it is not possible to choose $\lambda_{2}$ such
that $\partial_{x}\lambda_{2}>0$ for all $x$.  Therefore by
choosing $\lambda_{0}>0$ sufficiently large and
$\varepsilon,\delta>0$ sufficiently small, the matrix
$(\mathcal{B}_{ij})$ is positive definite.\par
  Although (2.3) yields the nondegeneracy conditions for $L_{t}$ in
(3.11) and (3.12), which then yield the above uniqueness proof, we
would like to show how (2.3) can be used directly to obtain this
goal.  According to (2.3) there exists a smooth function
$f(x,t)<0$ (assuming, as we may, that $M(x,t)<0$) such that
\begin{equation*}
\int_{0}^{l}(M^{-1}L_{t}(x,t)+f(x,t))dx=0,
\end{equation*}
where as in the proof of Lemma 3.2 the $x$-coordinate lies in the
range $[0,l)$.  Therefore the following equation admits a smooth
solution $\lambda_{2}$:
\begin{equation*}
2M\partial_{x}\lambda_{2}-\frac{3}{2}L_{t}=\frac{3}{2}Mf.
\end{equation*}
Then by choosing
\begin{equation*}
\lambda:=\lambda_{1}+\lambda_{2},\text{ }\text{ }\text{ }\text{
}\lambda_{1}=\frac{1}{4}\log\det I,
\end{equation*}
taking $\lambda_{0}>0$ sufficiently large, and $\delta>0$
sufficiently small, the calculations (3.16) and (3.17) show that
$(\mathcal{B}_{ij})$ is positive definite in $D_{\delta}$.\medskip

\textbf{Remark.}  \textit{Without (2.3) it is unclear if
$\mathcal{B}_{22}$ can be made positive.  Thus we may view (2.3)
as that which makes the system symmetric positive near a closed
asymptotic curve.}\medskip

  Lastly, we show that the boundary integral in (3.15) is nonnegative.
The boundary of $D_{\delta}$ consists of two curves, namely $t=0$
and $t=\delta$.  On $t=0$ the function $\overline{U}$ vanishes
identically, and on $t=\delta$ the unit normal has components
$\nu_{1}=0$, $\nu_{2}=1$.  Thus since $N>0$ and $L\leq 0$, the
desired conclusion follows.  Q.E.D.

\begin{center}
\textbf{4.  Proof of Theorems 1 and 2}
\end{center} \setcounter{equation}{0}
\setcounter{section}{4}

   We would now like to use Lemmas 3.1 and 3.3 to extend the
agreement of the second fundamental forms of $S$ and
$\overline{S}$ on $\partial\mathcal{C}$ to the whole of
$\mathcal{C}$.  Choose one family of asymptotic curves and
decompose $\mathcal{C}$ into the cylindrical domains
$\mathcal{C}_{1},\ldots,\mathcal{C}_{m}$ of Lemma 2.5.  We first
show that the second fundamental forms $II$ and $\overline{II}$
agree on $\mathcal{C}_{1}$, which contains the boundary curve
$\gamma_{1}$ of $\mathcal{C}$.  Let $\sigma_{1}$ be an asymptotic
curve of this family which emanates from $\gamma_{1}$ and spirals
toward $\partial\mathcal{C}_{1}$.  Let $\sigma_{2}$ be an
asymptotic curve of the other family which emanates from
$\gamma_{1}$, and let $p$ be a point at which they intersect.  Let
$R(p)$ denote the segment of $\sigma_{2}$ starting on $\gamma_{1}$
at $p_{0}$ and ending at $p$ on $\sigma_{1}$. Suppose that $II$
and $\overline{II}$ agree on $\gamma_{1}$.  By Lemma 3.1 they agree also
in a neighborhood of $p_{0}$.  Therefore, starting at points on
$R(p)$ near $p_{0}$, we may apply the Goursat uniqueness theorem
for the hyperbolic system (3.3), in which one prescribes data on
two intersecting characteristics and obtains a unique solution in
a small characteristic rectangle. This shows that $II$ and
$\overline{II}$ agree in a small characteristic rectangle which
contains a portion of $R(p)$.  Then working our way up $R(p)$ we
may repeatedly apply this procedure to show that $II$ and
$\overline{II}$ agree in a neighborhood of $R(p)$ all the way up
to $\sigma_{1}$.  Since this holds for any $p\in\sigma_{1}$ which
has an asymptotic curve of the other family connecting it to
$\gamma_{1}$, we conclude by the continuity method that $II$ and
$\overline{II}$ agree on all of $\sigma_{1}$ which lies below the
first closed asymptotic curve of the other family.  Since
$\sigma_{1}$ was arbitrary, we now have that the second
fundamental forms agree on all of the portion of $\mathcal{C}_{1}$
which is bounded by $\gamma_{1}$ and a closed asymptotic curve
$\Gamma$ of the other family.\par
   In order to continue this uniqueness beyond $\Gamma$, we apply
Lemma 3.3.  Note that $\Gamma$ cannot intersect
$\partial\mathcal{C}_{1}$, since if it does then an asymptotic
curve of the same family as $\Gamma$ must be tangent to
$\partial\mathcal{C}_{1}$ at some point, which is not possible. We
may now continue the method of the previous paragraph to obtain
uniqueness up until we hit another closed asymptotic curve of the
same family as $\Gamma$.  Eventually we will reach
$\partial\mathcal{C}_{1}$.  The same technique may be applied in
$\mathcal{C}_{2},\ldots,\mathcal{C}_{m}$.  It follows that
$II=\overline{II}$ on all of $\mathcal{C}$.

\begin{center}
\textbf{References}
\end{center}

\noindent 1.\hspace{.07in} A. D. Alexandrov, \textit{On a class of
closed surfaces}, Recuiel Math., $\mathbf{4}$ (1938),
69-77.\medskip

\noindent 2.\hspace{.07in} S. Cohn-Vossen, \textit{Unstarre
geschlossene Fl\"{a}chen}, Math. Annalen, $\mathbf{102}$
(1929-30), \par\hspace{-.01in} 10-29.\medskip

\noindent 3.\hspace{.07in} K. O. Friedrichs, \textit{Symmetric
positive linear differential equations}, Comm. Pure
\par\hspace{-.01in} Appl. Math., $\mathbf{1}$ (1958),
333-418.\medskip

\noindent 4.\hspace{.07in} N. Kuiper, \textit{Convex immersions of
closed surfaces in $E^{3}$}, Comment. Math. Hel-
\par\hspace{-.01in} vetici, $\mathbf{35}$ (1961), 85-92.\medskip

\noindent 5.\hspace{.07in} N. Kuiper, \textit{On surfaces in
Euclidean three-space}, Bull. Soc. Math. Belgique, $\mathbf{12}$
\par\hspace{-.01in} (1960), 5-22.\medskip

\noindent 6.\hspace{.07in} L. Nirenberg, \textit{Rigidity of a
class of closed surfaces}, 1963 Nonlinear Problems
\par\hspace{-.01in} (Proc. Sympos., Madison, Wis., 1962), Univ. of
Wisconsin Press, 177-193.\medskip

\noindent 7.\hspace{.07in} A. V. Pogorelov, \textit{Die eindeutige
Bestimmtheit allgemeiner konvexer Fl\"{a}chen} \par\hspace{-.01in}
(translation from the Russian), Akademie Verlag, Berlin,
1956.\medskip

\noindent 8.\hspace{.07in} M. Spivak, \textit{A Comprehensive
Introduction to Differential Geometry, Vol. 5},
\par\hspace{-.01in} Publish or Perish Inc., Houston, 1999.\medskip

\bigskip\bigskip\footnotesize

\noindent\textsc{Department of Mathematics, University of Notre
Dame, Notre Dame, IN 46556}\par

\noindent\textit{E-mail address}: \verb"Qing.Han.7@nd.edu"\bigskip

\noindent\textsc{Department of Mathematics, Stony Brook
University, Stony Brook, NY 11794}\par

\noindent\textit{E-mail address}: \verb"khuri@math.sunysb.edu"

\end{document}